\documentclass[12pt]{article}%
\usepackage{amsfonts}
\usepackage{sw20bams}
\usepackage{amsmath}
\usepackage{amssymb}
\usepackage{graphicx}%
\providecommand{\U}[1]{\protect\rule{.1in}{.1in}}
\begin{document}

\title{Small Matrices with Large Inverses: Unimodular $4\times4$ Cases}
\author{Steven Finch}
\date{August 17, 2026}
\maketitle

\begin{abstract}
How close to singularity can an $n\times n$ unimodular matrix be? \ For
ternary cases as $n$ increases, exact expressions are unlikely, but upon
fixing $n=4$ and assessing $(2k+1)$-ary cases as $k$ increases, we make
significant progress; similarly for $(k+1)$-ary cases of $4\times4$
nonnegative unimodular matrices.

\end{abstract}

\footnotetext{Copyright \copyright \ 2026 by Steven R. Finch. All rights
reserved.}

An $n\times n$ integer matrix $M$ with determinant $\pm1$ is called
\textbf{unimodular}. \ Let $\alpha=\left\Vert M\right\Vert $ and
$\beta=\left\Vert M^{-1}\right\Vert $, the maximum absolute entry of $M$ and
$M^{-1}$, respectively. \ If $n=2$, then $\alpha=\beta$ trivially. \ If
$n\geq2$, the ternary lower-triangular matrix%
\[
\left(
\begin{array}
[c]{cccccc}%
1 & 0 & 0 & \ldots & 0 & 0\\
-1 & 1 & 0 &  & 0 & 0\\
-1 & -1 & 1 &  & 0 & 0\\
\vdots &  &  & \ddots &  & \vdots\\
-1 & -1 & -1 &  & 1 & 0\\
-1 & -1 & -1 & \ldots & -1 & 1
\end{array}
\right)
\]
has $\alpha=1$, $\beta=2^{n-2}$. \ Define%
\[
\gamma_{n}^{-}(\alpha)=\max\left\{  \beta\geq1:\exists\text{ unimodular }%
M\in\mathbb{Z}^{n\times n}\text{ with }\left\Vert M\right\Vert =\alpha\text{,
}\left\Vert M^{-1}\right\Vert =\beta\right\}  .
\]
While $\gamma_{2}^{-}(k)=k$ for $k\geq1$ and $\gamma_{n}^{-}(1)=2^{n-2}$ for
$3\leq n\leq4$, other values $\gamma_{n}^{-}(k)$~are less clear. \ The
sequence $\gamma_{5}^{-}(1)=10$, $\gamma_{6}^{-}(1)=30$, $\gamma_{7}%
^{-}(1)=130$, ... defies description (we know only that $442\leq\gamma_{8}%
^{-}(1)\leq576$ \cite{EZ-matrx3, S1-matrx3}). \ Define also%
\[
\gamma_{n}^{+}(\alpha)=\max\left\{  \beta\geq1:\exists\text{ unimodular }%
M\in\mathbb{Z}_{\geq0}^{n\times n}\text{ with }\left\Vert M\right\Vert
=\alpha\text{, }\left\Vert M^{-1}\right\Vert =\beta\right\}  ;
\]
the alphabet underlying $M$ here has $k+1$ symbols rather than $2k+1$. \ The
sequence $\gamma_{3}^{+}(1)=1$, $\gamma_{4}^{+}(1)=2$, $\gamma_{5}^{+}(1)=3$,
$\gamma_{6}^{+}(1)=5$, $\gamma_{7}^{+}(1)=9$, $\gamma_{8}^{+}(1)=18$,
$\gamma_{9}^{+}(1)=42$, ... again defies description (we know $110\leq
\gamma_{10}^{+}(1)\leq144$ \cite{S2-matrx3}).

Our interest is in the case $n=4$. \ Let%
\[
\beta_{k}^{-}:=\gamma_{4}^{-}(k),\qquad\beta_{k}^{+}:=\gamma_{4}^{+}(k)
\]
and $A$ be a $3\times3$ submatrix of $M$ (delete one row and one column). For
$A\in\mathbb{Z}^{3\times3}$ let $g(A)$ be the greatest common divisor ($\gcd$)
of its nine $2\times2$ minors.

\textbf{Theorem. \ }Let $k\geq2$.

\textbf{(i) Nonnegative case.} If all entries of $M$ lie in $\{0,1,\dots,k\}$,
then
\[
|\det A|\ \leq\ \beta_{k}^{+}=k(2k^{2}-2k+1)=2k^{3}-2k^{2}+k.
\]

\textbf{(ii) Signed case. }If all entries of $M$ lie in $\{-k,\ldots
,-1,0,1,\ldots,k\}$, then
\[
|\det A|\ \leq\ \beta_{k}^{-}=k(2k-1)(2k+1)=4k^{3}-k.
\]

Both bounds are attained for every $k\geq2$.\medskip

Scenario (i) is \textit{not} a special case of (ii): a nonnegative matrix
obeys both bounds, but the sharp bound for the nonnegative subproblem is the
smaller $\beta_{k}^{+}$, with its own extremal family. The two are genuinely
distinct extremal problems sharing the scaffolding described below.

We write the rows of $A$ as $\mathbf{r}_{1},\mathbf{r}_{2},\mathbf{r}_{3}$,
and for a row index $i$ let $\mathbf{w}=\mathbf{w}^{(i)}$ be the
\textbf{cofactor vector} (the $2\times2$ minors involving the other two rows,
with signs), so $\det A=\mathbf{r}_{i}\cdot\mathbf{w}$. We assume WLOG\ that
$\det A>0$ throughout, since swapping two rows of $A$ negates the determinant
while leaving $\Vert A\Vert$ and every $|2\times2$ minor$|$ -- hence $g(A)$ -- unchanged.

\section{The common scaffolding}

\subsection{Multilinear cube bound}

The determinant is linear in each entry, so over the box its maximum is at a
vertex \cite{Rk-matrx3}, where it equals $k^{3}$ multiplied by the largest
$3\times3$ determinant with $0/1$ (respectively $0/\pm1$) entries
\cite{S1-matrx3, S2-matrx3}:
\[
\det A\leq2k^{3}\ \text{(nonnegative)},\qquad\det A\leq4k^{3}\ \text{(signed)}%
.
\]
Call these \textbf{Hadamard-type ceilings}. \ Likewise $\max_{\mathbf{r}_{i}%
}\mathbf{r}_{i}\cdot\mathbf{w}=k\Vert\mathbf{w}\Vert_{1}$ is itself such a
determinant, giving the \textbf{cofactor bounds}
\[
\Vert\mathbf{w}\Vert_{1}\leq2k^{2}\ \text{(nonnegative)},\quad\Vert
\mathbf{w}\Vert_{1}\leq4k^{2}\ \text{(signed)},\qquad|w_{j}|\leq2k^{2}.
\]

\subsection{Embedding}

This property is used identically in both \S 1 and \S 2.

\textbf{Lemma E.} If $M\in\mathbb{Z}^{n\times n}$ is unimodular and $A$ is any
$(n-1)\times(n-1)$ submatrix, then $g(A)=1$.

\textbf{Proof}. Deleting one row and one column is the same as conjugating by
permutation matrices (preserving unimodularity, $A$, and $g(A)$) and deleting
the last row/column. Hence assume
\[
M=%
\begin{pmatrix}
A & \mathbf{u}\\
\mathbf{v}^{\!\top} & d
\end{pmatrix}
,\qquad\det M=\pm1.
\]
The one-row/one-column \textbf{bordering identity} (a polynomial identity,
verified where $A$ is invertible via $A^{-1}=\operatorname{adj}(A)/\det A$) is%
\begin{equation}
\det M=d\det(A)-\mathbf{v}^{\!\top}\operatorname*{adj}(A)\mathbf{u}.
\tag{$\bigstar$}%
\end{equation}
The entries of $\operatorname{adj}(A)$ are $\pm$ the $2\times2$ minors, thus
$g(A)$ divides them all; and $\det A=\sum_{h}a_{1h}\operatorname{adj}(A)_{h1}$
shows $g(A)\mid\det A$. Therefore $g(A)\mid\det M=\pm1$, i.e., $g(A)=1$.
$\blacksquare$

Fixing the border data $(d,\mathbf{v},\mathbf{u})$ independently of $k$ in
$(\star)$ yields a \textbf{uniform family} $\{A_{k}\}$: one construction
embeds every $A_{k}$ as a $3\times3$ submatrix of a unimodular $M_{k}$ with
$\Vert M_{k}\Vert=k$ for all $k$ (cf. \S 1.4, \S 2.4).

\subsection{Shared architecture}

Each section is framed by:

\begin{itemize}
\item (Achievability) a statement about \textit{all} box matrices (no
unimodularity) showing the attainable values above the bound form a finite
"rigid" set;

\item (Embedding) each rigid survivor has $p\mid g(A)$ for some prime $p\mid
k$, contradicting Lemma E.
\end{itemize}

\section{Nonnegative case, $A\in\{0,\dots,k\}^{3\times3}$}

\subsection{Cross product}

$\qquad$

\textbf{Lemma C.} Given $\mathbf{p},\mathbf{q}\in\{0,\dots,k\}^{3}$,
$\mathbf{c}=\mathbf{p}\times\mathbf{q}$; if $P:=\sum_{i}\max(c_{i},0)\geq
T:=2k^{2}-2k+2$, then $k\mid c_{i}$ for all $i$.

\subsubsection{Proof}

Write the three signed minors explicitly:
\[
c_{1}=p_{2}q_{3}-p_{3}q_{2},\qquad c_{2}=p_{3}q_{1}-p_{1}q_{3},\qquad
c_{3}=p_{1}q_{2}-p_{2}q_{1}.
\]
Let $S=\{i:c_{i}>0\}$ and let $\mathbf{e}_{S}=\sum_{i\in S}\mathbf{e}_{i}%
\in\{0,1\}^{3}$ be its indicator vector. Since $P$ keeps exactly the positive
coordinates,
\begin{equation}
P=\sum_{i\in S}c_{i}=\mathbf{e}_{S}\cdot\mathbf{c}=\mathbf{e}_{S}%
\cdot(\mathbf{p}\times\mathbf{q})=\det(\mathbf{e}_{S},\mathbf{p},\mathbf{q}),
\tag{$\sharp$}%
\end{equation}
the last equality being the scalar-triple-product form of the determinant with
rows $\mathbf{e}_{S},\mathbf{p},\mathbf{q}$. This identity is what lets us
bound $P$ geometrically in the two easy cases.

The proof splits on $|S|\in\{0,1,2,3\}$. We show $|S|\leq1$ and $|S|=3$ are
\textit{vacuous} (they cannot reach the threshold $T$), and then analyze
$|S|=2$ exhaustively.

\subsubsection{Case $|S|\leq1$ (vacuous)}

Here $P$ is either $0$ (if $S=\varnothing$) or a single minor $c_{i}$. Any
single minor $c_{i}=p_{j}q_{\ell}-p_{\ell}q_{j}$ is a $2\times2$ determinant
with entries in $[0,k]$; its positive part is at most $p_{j}q_{\ell}\leq
k\cdot k=k^{2}$. Hence
\[
P\leq k^{2}.
\]
But $k^{2}<T$ because $T-k^{2}=k^{2}-2k+2=(k-1)^{2}+1>0$. So the hypothesis
$P\geq T$ is never met --\ nothing to prove.

\subsubsection{Case $|S|=3$ (vacuous)}

Now every $c_{i}>0$, so $\mathbf{e}_{S}=\mathbf{1}=(1,1,1)$ and by $(\sharp)$,
$P=\det(\mathbf{1},\mathbf{p},\mathbf{q})$. This determinant is multilinear in
the pair $(\mathbf{p},\mathbf{q})$ --\ linear in $\mathbf{p}$ for fixed
$\mathbf{q}$, and vice versa --\ so over the box $[0,k]^{3}\times
\lbrack0,k]^{3}$ its maximum is attained at a vertex, i.e., at $\mathbf{p}%
=k\mathbf{u},\ \mathbf{q}=k\mathbf{v}$ with $\mathbf{u},\mathbf{v}%
\in\{0,1\}^{3}$. Pulling the scalars out of two rows,
\[
P\leq\max_{\mathbf{u},\mathbf{v}\in\{0,1\}^{3}}\det(\mathbf{1},k\mathbf{u}%
,k\mathbf{v})=k^{2}\max_{\mathbf{u},\mathbf{v}\in\{0,1\}^{3}}\det
(\mathbf{1},\mathbf{u},\mathbf{v}).
\]
It remains to see that the $0/1$ determinant with a fixed all-ones row lies in
$\{-1,0,1\}$. Subtracting column $1$ from columns $2,3$ turns the top row
$(1,1,1)$ into $(1,0,0)$ and, expanding along it,
\[
\det(\mathbf{1},\mathbf{u},\mathbf{v})=\det\!%
\begin{pmatrix}
u_{2}-u_{1} & u_{3}-u_{1}\\
v_{2}-v_{1} & v_{3}-v_{1}%
\end{pmatrix}
.
\]
For a $0/1$ vector $\mathbf{u}$, the row $(u_{2}-u_{1},u_{3}-u_{1})$ lies in
$\{0,1\}^{2}$ if $u_{1}=0$ and in $\{-1,0\}^{2}$ if $u_{1}=1$ --\ it can never
mix $+1$ and $-1$. Checking the three sign-patterns (both rows nonnegative,
both nonpositive, or mixed) shows the $2\times2$ determinant is always in
$\{-1,0,1\}$. Hence $P\leq k^{2}<T$: again vacuous.

(Both vacuous cases confirm the intuition that spreading positive mass over
one or three coordinates cannot beat the two-coordinate "corner" configurations.)

\subsubsection{Case $|S|=2$ (the substantive case)}

Relabel so that $c_{1},c_{2}>0\geq c_{3}$; then $\mathbf{e}_{S}=(1,1,0)$ and,
by $(\sharp)$,
\[
P=c_{1}+c_{2}=\det\!%
\begin{pmatrix}
1 & 1 & 0\\
p_{1} & p_{2} & p_{3}\\
q_{1} & q_{2} & q_{3}%
\end{pmatrix}
.
\]
Expanding along the top row,
\[
P=(p_{2}q_{3}-p_{3}q_{2})+(p_{3}q_{1}-p_{1}q_{3})=q_{3}\underset{_{=:\,\xi}%
}{\underbrace{(p_{2}-p_{1})}}+p_{3}\underset{_{=:\,\eta}}{\underbrace
{(q_{1}-q_{2})}}=q_{3}\xi+p_{3}\eta.
\]

\paragraph{Each product is capped at $k^{2}$}

We have $q_{3},p_{3}\in\{0,\dots,k\}$ while $\xi,\eta\in\{-k,\dots,k\}$. Since
$P\geq T>0$, the two products cannot both be nonpositive; in fact we show each
is \textit{large}. Each obeys $q_{3}\xi\leq k\cdot k=k^{2}$ and $p_{3}\eta\leq
k^{2}$. If one of them, say $q_{3}\xi$, were $\leq(k-1)^{2}$, then
\[
P\leq(k-1)^{2}+k^{2}=2k^{2}-2k+1=T-1<T,
\]
contradicting the hypothesis. Hence \textit{both} products satisfy
\begin{equation}
q_{3}\xi\ \geq\ (k-1)^{2}+1,\qquad p_{3}\eta\ \geq\ (k-1)^{2}+1. \tag{$\dag$}%
\end{equation}

\paragraph{Factor-pairs which achieve $(\dag)$}

In particular each product is positive, so $\xi,\eta>0$; thus $q_{3},\xi
,p_{3},\eta$ are all integers in $\{1,\dots,k\}$. We need two such integers
with product $\geq(k-1)^{2}+1=k^{2}-2k+2$. Testing the possibilities (both
factors $\leq k$):

\begin{itemize}
\item $k\cdot k=k^{2}$ $\ \checkmark$;

\item $k(k-1)=k^{2}-k\geq k^{2}-2k+2\iff k\geq2$ $\ \checkmark$;

\item $k(k-2)=k^{2}-2k<k^{2}-2k+2$ $\ \times$;

\item $(k-1)^{2}=k^{2}-2k+1<k^{2}-2k+2$ \ $\times$.
\end{itemize}

\noindent So the only admissible pairs are
\[
(q_{3},\xi),\ (p_{3},\eta)\ \in\ \{(k,k),\,(k,k-1),\,(k-1,k)\},
\]
with product $k^{2}$ for $(k,k)$ and $k^{2}-k$ for the other two.

\paragraph{The sum couples the two pairs}

The total is $P=q_{3}\xi+p_{3}\eta$, a sum of two values each in
$\{k^{2},\,k^{2}-k\}$, hence $P\in\{2k^{2},\,2k^{2}-k,\,2k^{2}-2k\}$. The
threshold $P\geq T=2k^{2}-2k+2$ excludes the smallest option $2k^{2}-2k$
--\ so \textit{at least one} of the two pairs must be $(k,k)$ (contributing
the full $k^{2}$). This is exactly the admissibility table:%

\[%
\begin{tabular}
[c]{|l|l|}\hline
$(q_{3},\xi)$ & admissible $(p_{3},\eta)$\\\hline
$(k,k)$ & $(k,k),(k,k-1),(k-1,k)$\\\hline
$(k,k-1)$ & $(k,k)$ only\\\hline
$(k-1,k)$ & $(k,k)$ only\\\hline
\end{tabular}
\
\]

\paragraph{Decoding $(q_{3},\xi,p_{3},\eta)$ back to $(\mathbf{p},\mathbf{q}%
)$}

We now recover the actual rows. The values $q_{3},p_{3}$ are read off
directly. For $\xi=p_{2}-p_{1}$ with $p_{1},p_{2}\in\{0,\dots,k\}$:
\[
\xi=k\ \Rightarrow\ (p_{1},p_{2})=(0,k)\ \text{(unique)},
\]%
\[
\xi=k-1\ \Rightarrow\ (p_{1},p_{2})\in\{(0,k-1),(1,k)\}\ \text{(two choices)}%
\]
and identically for $\eta=q_{1}-q_{2}$. So each "$k$" entry decodes uniquely,
while each "$k-1$" entry branches in two. Counting decodings across the
admissibility table --\ $(k,k)\&(k,k)$ gives $1$; $(k,k)\&(k,k-1)$ gives $2$;
$(k,k)\&(k-1,k)$ gives $1$; $(k,k-1)\&(k,k)$ gives $2$; $(k-1,k)\&(k,k)$ gives
$1$ --\ produces exactly $1+2+1+2+1=7$ configurations, which are the seven
rows below (with $c_{1}=p_{2}q_{3}-p_{3}q_{2}$, etc.):%

\[%
\begin{tabular}
[c]{|l|l|l|l|l|}\hline
\# & $(p_{1},p_{2},p_{3})$ & $(q_{1},q_{2},q_{3})$ & $(c_{1},c_{2},c_{3})$ &
$P$\\\hline
1 & $(0,k,k)$ & $(k,0,k)$ & $(k^{2},k^{2},-k^{2})$ & $2k^{2}$\\\hline
2 & $(0,k,k)$ & $(k-1,0,k)$ & $(k^{2},k^{2}-k,-(k^{2}-k))$ & $2k^{2}%
-k$\\\hline
3 & $(0,k,k)$ & $(k,1,k)$ & $(k^{2}-k,k^{2},-k^{2})$ & $2k^{2}-k$\\\hline
4 & $(0,k,k-1)$ & $(k,0,k)$ & $(k^{2},k^{2}-k,-k^{2})$ & $2k^{2}-k$\\\hline
5 & $(0,k-1,k)$ & $(k,0,k)$ & $(k^{2}-k,k^{2},-(k^{2}-k))$ & $2k^{2}%
-k$\\\hline
6 & $(1,k,k)$ & $(k,0,k)$ & $(k^{2},k^{2}-k,-k^{2})$ & $2k^{2}-k$\\\hline
7 & $(0,k,k)$ & $(k,0,k-1)$ & $(k^{2}-k,k^{2},-k^{2})$ & $2k^{2}-k$\\\hline
\end{tabular}
\
\]

For instance, branch 1 comes from $(q_{3},\xi)=(p_{3},\eta)=(k,k)$: $\xi=k$
gives $(p_{1},p_{2})=(0,k)$, $\eta=k$ gives $(q_{1},q_{2})=(k,0)$, and then
$c_{1}=k\cdot k-k\cdot0=k^{2}$, $c_{2}=k\cdot k-0=k^{2}$, $c_{3}=0-k\cdot
k=-k^{2}$. Branches 2,3 are the two decodings of $(k,k)\&(k,k-1)$, branches
5,6 the two decodings of $(k,k-1)\&(k,k)$, and 4,7 the unique decodings
involving a $(k-1,k)$ pair.

\paragraph{Conclusion of the case}

In every one of the seven branches each coordinate $c_{i}$ takes a value in
$\{k^{2},\ k^{2}-k,\ -(k^{2}-k),\ -k^{2}\}$, and each of these is a multiple
of $k$. Hence $k\mid c_{1},c_{2},c_{3}$.

Combining all cases, whenever $P\ge T$ we are necessarily in the $|S|=2$
situation, where $k\mid c_{i}$ for all $i$. $\blacksquare$

\subsubsection{Sharpness}

The threshold $T$ cannot be lowered. Take $\mathbf{p}=(0,k,k-1)$,
$\mathbf{q}=(k,1,k)$: then $\xi=k,\ p_{3}=k-1,\ \eta=k-1,\ q_{3}=k$, so
\[
P=q_{3}\xi+p_{3}\eta=k^{2}+(k-1)^{2}=2k^{2}-2k+1=T-1,
\]
one below threshold, and
\[
c_{1}=p_{2}q_{3}-p_{3}q_{2}=k\cdot k-(k-1)\cdot1=k^{2}-k+1\equiv
1\,\,(\operatorname{mod}k),
\]
so $k\nmid c_{1}$. Thus $P\geq T$ (not merely $P\geq T-1$) is exactly what the
divisibility conclusion requires.

\subsection{Common factors}

$\qquad$

\textbf{Minor Rigidity}

If $\det A>\beta_{k}^{+}$, then $k\mid$ every $2\times2$ minor, i.e., $k\mid
g(A)$.

\textbf{Proof}. For an ordered pair of rows the cross product
\[
\mathbf{r}_{j}\times\mathbf{r}_{\ell}=(m_{1},m_{2},m_{3})
\]
collects the three signed $2\times2$ minors built from rows $j,\ell$ (one per
omitted column). Write
\[
P^{(i)}:=\sum_{h=1}^{3}\max\left\{  m_{h},0\right\}
\]
for the \textbf{positive-part sum} of the minors of the complementary pair
$\{j,\ell\}$, taken in an orientation that makes the triple product equal to
$\det A$.

Cofactor (row-$i$) expansion is exactly the scalar triple product. Choosing
the cyclically ordered complementary pair:
\[
i=1\!:(j,\ell)=(2,3),\quad i=2\!:(j,\ell)=(3,1),\quad i=3\!:(j,\ell)=(1,2),
\]
all of which are \textit{even} permutations of $(1,2,3)$, we obtain for each
$i$
\[
\det A=\mathbf{r}_{i}\cdot(\mathbf{r}_{j}\times\mathbf{r}_{\ell})=\sum
_{h=1}^{3}(r_{i})_{h}\,m_{h}.
\]
Now $(r_{i})_{h}\in\lbrack0,k]$, so each term obeys $(r_{i})_{h}$ $m_{h}\leq
k\max(m_{h},0)$: when $m_{h}>0$ the factor is at most $k$, and when $m_{h}%
\leq0$ the term is $\leq0$ (note that $(r_{i})_{h}\geq0$ is essential here).
Summing,
\[
\ \det A\ \leq\ k\,P^{(i)}\quad\text{for each }i=1,2,3.\
\]

Assume $\det A>\beta_{k}^{+}$. Since $\det A\in\mathbb{Z}$,
\[
\det A\ \geq\ \beta_{k}^{+}+1=k(2k^{2}-2k+1)+1.
\]
Combine with the preceding, for \textit{each} $i$,
\[
P^{(i)}\ \geq\ \frac{\det A}{k}\ \geq\ \frac{k(2k^{2}-2k+1)+1}{k}%
=2k^{2}-2k+1+\dfrac{1}{k}\ >\ 2k^{2}-2k+1.
\]
But $P^{(i)}$ is an integer, so the strict inequality lifts to
\[
P^{(i)}\ \geq\ 2k^{2}-2k+2=T\qquad(i=1,2,3).
\]

For a fixed $i$, $P^{(i)}\geq T$ is precisely the hypothesis of Lemma C for
the row pair $\{j,\ell\}$. Hence $k$ divides all three minors of $\{j,\ell\}$.
\ The three choices of $i$ produce the three distinct row pairs
$\{2,3\},\{3,1\},\{1,2\}$, and each pair contributes its three minors (one per
omitted column). That is
\[
3\ \text{pairs}\times3\ \text{minors}=9=\binom{3}{2}\binom{3}{2}%
\]
minors in total --\ i.e., \textit{every} $2\times2$ minor of $A$. Therefore
$k$ divides each one, and a fortiori
\[
k\ \mid\ \gcd(\text{all }2\times2\text{ minors})=g(A).\qquad\blacksquare
\]

\subsection{Divisibility upgrade}

$\qquad$

\textbf{Spectral Rigidity}

$\det A>\beta_{k}^{+}$ implies $k^{2}\mid\det A$. Consequently the only
attainable determinant values in the window $(\beta_{k}^{+},2k^{3}]$ are
$2k^{3}-k^{2}$ and $2k^{3}$.

\textbf{Proof} (simultaneous strong induction on $k$). The conclusion is
really \textit{two} statements bundled together:

\begin{itemize}
\item (D) the divisibility claim $\;k^{2}\mid\det A$;

\item (R) the rigidity claim that the only attainable values in $(\beta
_{k}^{+},2k^{3}]$ are $2k^{3}-k^{2}$ and $2k^{3}$.
\end{itemize}

\noindent These statements are not proved one-after-the-other; they are
carried \textit{together }through the induction, each feeding the other:%

\[
\underset{\text{induction hypothesis}}{\underbrace{(R)\text{ at }k_{1}<k}%
}\ \longrightarrow\ \underset{_{\det A=2k^{3}\text{ forced}}}{\underbrace
{(D)\text{ at }k}\ }\longrightarrow\ (R)\text{ at }k.
\]

Let $s_{1}\mid s_{2}\mid s_{3}$ be the Smith invariant factors of $A$. Two
standard facts are:

\begin{itemize}
\item $\det A=s_{1}s_{2}s_{3}$ (since $\det A>0$, there is no sign issue);

\item the gcd of the $2\times2$ minors equals the product of the first two
invariant factors, $g(A)=s_{1}s_{2}.$
\end{itemize}

Minor Rigidity applies because $\det A>\beta_{k}^{+}$; it yields
\[
k\mid g(A)=s_{1}s_{2}.
\]

For a nonzero integer $n$ and a fixed prime $p$, the $p$\textbf{-adic
valuation} $v_{p}(n)$ is the exponent of the highest power of $p$ dividing
$n$:
\[
v_{p}(n)=\max\{\,e\geq0:p^{e}\mid n\,\}.
\]
By convention $v_{p}(0)=+\infty$. \ The two properties we shall employ are:
\[
v_{p}(mn)=v_{p}(m)+v_{p}(n),\qquad p^{e}\mid n\iff v_{p}(n)\geq e.
\]

We must show that $k^{2}\mid\det A=s_{1}s_{2}s_{3}$. \ It suffices to prove,
prime by prime, that
\[
v_{p}(\det A)\ \geq\ 2\,v_{p}(k)\qquad\text{for every prime }p\mid k,
\]
because then $v_{p}(\det A)\geq v_{p}(k^{2})$ for all $p$, i.e., $k^{2}%
\mid\det A$.

Fix $p\mid k$ and write
\[
\kappa:=v_{p}(k),\qquad a:=v_{p}(s_{1})\leq b:=v_{p}(s_{2})\leq c:=v_{p}%
(s_{3}),
\]
the inequalities holding because $s_{1}\mid s_{2}\mid s_{3}$. The hypothesis
$k\mid s_{1}s_{2}$ reads
\begin{equation}
a+b\ \geq\ \kappa. \tag{1}%
\end{equation}
Our goal for this prime is $a+b+c\geq2\kappa$. Suppose for contradiction that
\begin{equation}
a+b+c<2\kappa. \tag{2}%
\end{equation}
Since $c\geq b$,
\begin{equation}
a+2b\ \leq\ a+b+c\ <\ 2\kappa. \tag{3}%
\end{equation}
by (2). Subtracting (1) from (3):
\begin{equation}
b=(a+2b)-(a+b)<2\kappa-\kappa=\kappa\tag{4}%
\end{equation}
so $b<\kappa$. \ If we had $a=0$, then (1) would give $b\geq\kappa$,
contradicting (4). Hence
\begin{equation}
a\geq1,\quad\text{i.e., }p\mid s_{1}. \tag{5}%
\end{equation}

The first invariant factor $s_{1}$ is exactly the gcd of all entries of $A$.
By (5), $p$ divides every entry, so
\begin{equation}
A=p\,A^{\prime},\qquad A^{\prime}\in\{0,1,\dots,k_{1}\}^{3\times3},\quad
k_{1}:=k/p<k, \tag{6}%
\end{equation}
and $A^{\prime}$ is again a nonnegative integer matrix, of the required form
for parameter $k_{1}$.

From $A=pA^{\prime}$,
\begin{equation}
\det A=p^{3}\det A^{\prime},\qquad\text{so}\qquad\det A^{\prime}=\frac{\det
A}{p^{3}}>\frac{\beta_{k}^{+}}{p^{3}}. \tag{7}%
\end{equation}
Compute the right-hand side with $k=pk_{1}$:
\[
\frac{\beta_{k}^{+}}{p^{3}}=\frac{2k^{3}}{p^{3}}-\frac{2k^{2}}{p^{3}}+\frac
{k}{p^{3}}=2k_{1}^{3}-\frac{2k_{1}^{2}}{p}+\frac{k_{1}}{p^{2}}.
\]
Now compare with $2k_{1}^{3}-k_{1}^{2}$. Because $p\geq2$ we have
$\tfrac{2k_{1}^{2}}{p}\leq k_{1}^{2}$, hence%
\begin{equation}
\frac{\beta_{k}^{+}}{p^{3}}=2k_{1}^{3}-\frac{2k_{1}^{2}}{p}+\frac{k_{1}}%
{p^{2}}\ \geq\ 2k_{1}^{3}-k_{1}^{2}+\frac{k_{1}}{p^{2}}\ >\ 2k_{1}^{3}%
-k_{1}^{2}. \tag{8}%
\end{equation}
(The middle term $\tfrac{k_{1}}{p^{2}}>0$ makes the last inequality strict,
even when $p=2$ makes the first an equality.)

Combining (7) and (8):
\begin{equation}
\ \det A^{\prime}>2k_{1}^{3}-k_{1}^{2}.\ \tag{9}%
\end{equation}
Also, $2k_{1}^{3}-k_{1}^{2}>\beta_{k_{1}}^{+}=2k_{1}^{3}-2k_{1}^{2}+k_{1}$
exactly when $k_{1}^{2}>k_{1}$, i.e., $k_{1}\geq2$; and for $k_{1}=1$ one
checks directly that $2k_{1}^{3}-k_{1}^{2}=1=\beta_{1}^{+}$, while (9) already
gives $\det A^{\prime}>1$. In every case $\det A^{\prime}$ lands
\textit{strictly above} $2k_{1}^{3}-k_{1}^{2}$ and inside the domain of the
induction hypothesis.

By (7)-(9) we have $\det A^{\prime}>\beta_{k_{1}}^{+}$, so the induction
hypothesis (Spectral Rigidity at $k_{1}$) tells us the \textit{only} possible
values in $(\beta_{k_{1}}^{+},\,2k_{1}^{3}]$ are $2k_{1}^{3}-k_{1}^{2}$ and
$2k_{1}^{3}$. But (9) excludes the smaller one, and $\det A^{\prime}\leq
2k_{1}^{3}$ by the Hadamard-type ceiling. Therefore
\begin{equation}
\det A^{\prime}=2k_{1}^{3}. \tag{10}%
\end{equation}

From (6) and (10),
\[
\det A=p^{3}\det A^{\prime3}\,(2k_{1}^{3})=2(pk_{1})^{3}=2k^{3}.
\]
Hence
\[
a+b+c=v_{p}(\det A)=v_{p}(2k^{3})=v_{p}(2)+3\kappa\ \geq\ 3\kappa
\ \geq\ 2\kappa,
\]
the last step using $\kappa=v_{p}(k)\geq1$. This contradicts the supposition (2).

Thus $a+b+c\geq2\kappa$ for the chosen prime, and since $p\mid k$ was
arbitrary,
\[
v_{p}(\det A)\geq2\,v_{p}(k)\ \text{ for all }p\mid k
\]
and therefore $k^{2}\mid\det A.$

We finally must pin down the two surviving values, equipped with knowledge
that $\det A$ is a multiple of $k^{2}$ lying in $(\beta_{k}^{+},\,2k^{3}]$.
Write $\det A=jk^{2}$. Then
\[
\beta_{k}^{+}<jk^{2}\leq2k^{3}\;\Longleftrightarrow\;2k-2+\tfrac{1}%
{k}\;<\;j\;\leq\;2k.
\]
Since $2k-2<2k-2+\tfrac{1}{k}<2k-1$, the integer $j$ must satisfy
$j\in\{2k-1,\,2k\}$, giving
\begin{equation}
\det A\in\{(2k-1)k^{2},\,2k\cdot k^{2}\}=\{\,2k^{3}-k^{2},\ 2k^{3}\,\}.
\tag{$\natural$}%
\end{equation}

For $k=1$, $A$ is a $0/1$ matrix; the maximum $3\times3$ determinant is $2$,
and $\beta_{1}^{+}=1$. The window $(1,2]$ contains only the value
$2=2\cdot1^{3}$, and $k^{2}=1$ divides everything -- so the statement holds
and anchors the induction (every prime reduction in (7)-(9) eventually
descends to $k_{1}=1$).

For $k=2$, we have $\beta_{2}^{+}=2\cdot8-2\cdot4+2=10$, window $(10,16]$, and
the multiples of $k^{2}=4$ there are exactly $12=2k^{3}-k^{2}$ and $16=2k^{3}%
$, matching ($\natural$). $\blacksquare$

\subsection{\textbf{Bound and tightness (}$+$)}

If $|\det A|>\beta_{k}^{+}$ then $\det A\in\{2k^{3}-k^{2},2k^{3}\}$, both with
$k\mid g(A)$ by Minor Rigidity, contradicting Lemma E. Thus $|\det A|\leq
\beta_{k}^{+}$.

Let $G$ denote the group of transformations $A\mapsto P\,A\,Q$ with $P$, $Q$
permutation matrices; $G$ preserves nonnegativity, $\Vert A\Vert$, $|\det A|$
and $g(A)$. The two $G$-orbit representatives
\[
A_{k}^{+}=%
\begin{pmatrix}
0 & k & k-1\\
k & 1 & k\\
k & k & 0
\end{pmatrix}
,\qquad B_{k}^{+}=%
\begin{pmatrix}
0 & k & k-1\\
k & 0 & k\\
k & k-1 & 0
\end{pmatrix}
\]
have $\det=\beta_{k}^{+}$, $g=1$, with adjugate first rows $(-k^{2}%
,k^{2}-k,k^{2}-k+1)$ and $(-k^{2}+k,(k-1)^{2},k^{2})$. Taking
$d=0,\ \mathbf{v}=(1,0,0)$ and $\mathbf{u}_{A}=(0,-1,1),\ \mathbf{u}%
_{B}=(2,1,1)$ in $(\star)$ gives $\mathbf{v}^{\!\top}\operatorname{adj}%
(A_{k}^{+})\mathbf{u}_{A}=1=\mathbf{v}^{\!\top}\operatorname{adj}(B_{k}%
^{+})\mathbf{u}_{B}$, $\det(M_{A}^{+})=-1=\det(M_{B}^{+})$:
\[
M_{A}^{+}=%
\begin{pmatrix}
0 & k & k-1 & 0\\
k & 1 & k & -1\\
k & k & 0 & 1\\
1 & 0 & 0 & 0
\end{pmatrix}
,\quad M_{B}^{+}=%
\begin{pmatrix}
0 & k & k-1 & 2\\
k & 0 & k & 1\\
k & k-1 & 0 & 1\\
1 & 0 & 0 & 0
\end{pmatrix}
,\quad\Vert M\Vert=k.
\]

\section{Signed case, $A\in\{-k,\dots,k\}^{3\times3}$}

\subsection{Minor gaps}

$\qquad$

\textbf{Lemma P.} No integer in $(2k^{2}-k,\ 2k^{2})$ is a $2\times2$ minor of
a matrix with entries in $[-k,k]$.

\textbf{Proof.} A minor is $P-Q$, $P=\alpha\delta,\ Q=\beta\gamma\in
\lbrack-k^{2},k^{2}]$. If $P-Q>2k^{2}-k$ then $(k^{2}-P)+(k^{2}+Q)<k$, a sum
of two nonnegative integers. The largest product $\leq k^{2}$ is $k^{2}$; the
next is $k^{2}-k$. Hence each of $k^{2}-P,\ k^{2}+Q$ is $0$ or $\geq k$; their
sum $<k$ forces both $0$, i.e., $P-Q=2k^{2}$. $\blacksquare\medskip$

Therefore every attainable $|w_{j}|$ lies in $\{0,\dots,2k^{2}-k\}\cup
\{2k^{2}\}$.

\subsection{Achievability}

$\qquad$

\textbf{Lemma B.} The band $\{4k^{3}-k+1,\dots,4k^{3}-1\}$ is empty.

\textbf{Proof.} Fix the expansion row -- say row $1$, with entries
$a_{1},a_{2},a_{3}$ -- and let
\[
\mathbf{w}=(w_{1},w_{2},w_{3})=\mathbf{r}_{2}\times\mathbf{r}_{3}%
\]
be the vector of signed $2\times2$ minors. Then cofactor expansion is the
triple product
\[
\det A=\sum_{j=1}^{3}a_{j}\,w_{j}.
\]

We use three facts established earlier:

\begin{itemize}
\item (N) the \textit{minor-norm ceiling}: any realizable minor vector has
$\Vert\mathbf{w}\Vert_{1}\leq4k^{2}$;

\item (M) the \textit{single-minor bound}: each signed $2\times2$ minor
satisfies $|w_{j}|\leq2k^{2}$;

\item (P) the $2\times2$\textit{ minor gap}: $\{2k^{2}-k+1,\dots,2k^{2}-1\}$
contains no attainable signed minor; i.e., an attainable minor with
$|w_{j}|\geq2k^{2}-k+1$ must equal $2k^{2}$.\medskip
\end{itemize}

For each $j$ with $w_{j}\neq0$ put
\[
c_{j}:=k-a_{j}\,\operatorname{sgn}(w_{j}).
\]
Since $|a_{j}|\leq k$ we have $a_{j}\operatorname{sgn}(w_{j})\leq k$, hence
$c_{j}\geq0$; it is an integer. Define the \textbf{deficit}
\[
m:=\sum_{j}c_{j}|w_{j}|.
\]
Using $\operatorname{sgn}(w_{j})\,|w_{j}|=w_{j}$,
\begin{equation}
m=\sum_{j}\left(  k-a_{j}\operatorname{sgn}(w_{j})\right)  |w_{j}|=k\sum
_{j}|w_{j}|-\sum_{j}a_{j}w_{j}=k\Vert\mathbf{w}\Vert_{1}-\det A. \tag{$\flat$}%
\end{equation}
So $m\geq0$ is exactly the gap between $\det A$ and its row-expansion ceiling
$k\Vert\mathbf{w}\Vert_{1}$.

Now suppose, for contradiction,
\[
\det A=4k^{3}-t,\qquad1\leq t\leq k-1.
\]
From $(\flat)$ with $m\geq0$,
\[
k\Vert\mathbf{w}\Vert_{1}\geq\det A=4k^{3}-t\quad\Longrightarrow\quad
\Vert\mathbf{w}\Vert_{1}\geq4k^{2}-\tfrac{t}{k}.
\]
Because $1\leq t\leq k-1$ we have $0<\tfrac{t}{k}<1$, so $\Vert\mathbf{w}%
\Vert_{1}>4k^{2}-1$. But $\Vert\mathbf{w}\Vert_{1}$ is a nonnegative integer
and, by (N), is $\leq4k^{2}$. The only integer in $(4k^{2}-1,\,4k^{2}]$ is
$4k^{2}$, so
\begin{equation}
\Vert\mathbf{w}\Vert_{1}=4k^{2},\qquad\text{whence by }(\flat)\quad
m=k\cdot4k^{2}-(4k^{3}-t)=t\leq k-1. \tag{1'}%
\end{equation}

\subsubsection{Aligned vs. small coordinates}

Classify each coordinate with $w_{j}\ne0$:

\begin{itemize}
\item \textbf{aligned}: $c_{j}=0$, i.e., $a_{j}=k\operatorname{sgn}(w_{j})$
(the expansion entry is $\pm k$, matching the minor's sign); contributes $0$
to $m$.

\item \textbf{small}: $c_{j}\geq1$. Then $m=\sum_{i}c_{i}|w_{i}|\geq
c_{j}|w_{j}|\geq|w_{j}|$, so by (1'), $|w_{j}|\leq m\leq k-1$.
\end{itemize}

\noindent Moreover, since every small coordinate has $c_{j}\geq1$,
\begin{equation}
\sum_{j\text{ small}}|w_{j}|\ \leq\ \sum_{j\text{ small}}c_{j}|w_{j}%
|\ \leq\ m\ \leq\ k-1, \tag{2'}%
\end{equation}
and therefore, using $\Vert\mathbf{w}\Vert_{1}=4k^{2}$,
\begin{equation}
\sum_{j\text{ aligned}}|w_{j}|=4k^{2}-\sum_{j\text{ small}}|w_{j}%
|\ \geq\ 4k^{2}-(k-1). \tag{3'}%
\end{equation}

\subsubsection{Exactly two aligned, one small}

$\qquad\qquad\qquad\qquad$

\textit{At least two aligned.} If at most one coordinate were aligned, then by
(M) the aligned sum would be $\leq2k^{2}$. But
\[
2k^{2}<4k^{2}-(k-1)\iff0<2k^{2}-k+1,
\]
which always holds, contradicting (3'). Hence $\geq2$ coordinates are aligned.

\textit{At least one small.} If every coordinate were aligned then $m=0$,
contradicting $m=t\geq1$.

With only three coordinates, "$\geq2$ aligned and $\geq1$ small" forces
exactly two aligned and one small, all three nonzero (there is no room left
for a zero coordinate). Write the small index as $c$; then $w_{c}\neq0$ (small
coordinates are nonzero by classification) and $|w_{c}|\leq m\leq k-1$ (by
preceding small bullet), so
\begin{equation}
1\leq|w_{c}|\leq k-1. \tag{4'}%
\end{equation}

\subsubsection{Contradiction via (P)}

Let $|w_{a}|,|w_{b}|$ be the two aligned values. By (3') and the exact count,
\begin{equation}
|w_{a}|+|w_{b}|=\Vert\mathbf{w}\Vert_{1}-|w_{c}|=4k^{2}-|w_{c}|. \tag{5'}%
\end{equation}
Each aligned value is $\leq2k^{2}$ by (M), so from (5'),
\[
|w_{a}|=(|w_{a}|+|w_{b}|)-|w_{b}|\ \geq\ (4k^{2}-|w_{c}|)-2k^{2}=2k^{2}%
-|w_{c}|\ \geq\ 2k^{2}-(k-1)=2k^{2}-k+1,
\]
and symmetrically $|w_{b}|\geq2k^{2}-k+1$. Thus both aligned minors lie in
$[\,2k^{2}-k+1,\ 2k^{2}\,]$. By the $2\times2$ gap (P), each must equal the
top value:
\[
|w_{a}|=|w_{b}|=2k^{2}.
\]
Substituting into (5'),
\[
|w_{c}|=4k^{2}-(|w_{a}|+|w_{b}|)=4k^{2}-4k^{2}=0,
\]
contradicting $|w_{c}|\geq1$ from (4').

No box matrix can have $\det A=4k^{3}-t$ with $1\leq t\leq k-1$; i.e., the
width-$(k-1)$ band $\{4k^{3}-k+1,\dots,4k^{3}-1\}$ is empty. $\blacksquare$

\subsection{Top value $4k^{3}$ is excluded}

$\qquad$

\textbf{Lemma V.} $\det A=4k^{3},\ \Vert A\Vert\leq k$ jointly imply that
every $2\times2$ minor $\equiv0$ $\left(  \operatorname{mod}p\right)  $ for
each prime $p\mid k$; i.e., $p\mid g(A)$.

\textbf{Proof.} Start with the Laplace expansion along row $i$:
\[
\det A=\sum_{j}a_{ij}w_{j}^{(i)},
\]
where $w_{j}^{(i)}$ are the row-$i$ cofactors (i.e., the signed $2\times2$
minors). We have the chain
\[
4k^{3}=\sum_{j}a_{ij}w_{j}^{(i)}\;\leq\;k\sum_{j}|w_{j}^{(i)}|=k\,\Vert
\mathbf{w}^{(i)}\Vert_{1}%
\]
giving $\Vert\mathbf{w}^{(i)}\Vert_{1}\geq4k^{2}$. To turn the chain into a
string of equalities, we also need an upper bound $\Vert\mathbf{w}^{(i)}%
\Vert_{1}\leq4k^{2}$. The three cofactors of a row are the components of
$u\times v$ for the other two rows $u,v$, and for fixed $v$ the map
$u\mapsto\Vert u\times v\Vert_{1}$ is convex, so its maximum over $\Vert
u\Vert_{\infty},\Vert v\Vert_{\infty}\leq k$ is attained at vertices
$u,v\in\{\pm k\}^{3}$; it follows that $\Vert u\times v\Vert_{1}\leq4k^{2}$.
From this,
\[
4k^{3}\leq k\Vert\mathbf{w}^{(i)}\Vert_{1}\leq4k^{3}%
\]
collapses, giving simultaneously $\Vert\mathbf{w}^{(i)}\Vert_{1}=4k^{2}$ and
equality in $\sum a_{ij}w_{j}\leq k\sum|w_{j}|$, which forces $a_{ij}=\pm k$
(sign aligned with $w_{j}^{(i)}$) wherever $w_{j}^{(i)}\neq0$. Since $p\mid
k$, those entries are $\equiv0$.

With $\Vert\mathbf{w}^{(i)}\Vert_{1}=4k^{2}$ and each $|w_{j}|\leq2k^{2}$, two
zeros would force the third to be $4k^{2}>2k^{2}$. Thus each row has $\geq2$
nonzero cofactors.

If row $i_{0}$ has its single zero cofactor, say, at column $2$, then
$|w_{1}^{(i_{0})}|+|w_{3}^{(i_{0})}|=4k^{2}$ with each $\leq2k^{2}$, so both
equal $2k^{2}$. A minor of value $\pm2k^{2}$ satisfies $|ad-bc|=2k^{2}%
=|ad|+|bc|$, forcing $|a|=|b|=|c|=|d|=k$, all four entries $\pm k$. The minor
$w_{1}^{(i_{0})}$ lives in columns $\{2,3\}$ and $w_{3}^{(i_{0})}$ in columns
$\{1,2\}$, both in the two rows other than $i_{0}$; their column sets cover
$\{1,2,3\}$ (any two distinct $2$-subsets of a $3$-set have full union --
which is why the choice of column $2$ is genuinely WLOG). Hence those two rows
are entirely $\pm k\equiv0$, only row $i_{0}$ survives mod $p$, so
$\operatorname{rank}\leq1$.

If no row has a zero cofactor, then all nine entries are $\pm k\equiv0$.
Either way we have $p\mid$ every minor. $\blacksquare\medskip$

The two cases are exhaustive, so in all cases every minor $\equiv0$ $\left(
\operatorname{mod}p\right)  $ , i.e., $p\mid g(A)$. Combined with Lemma E
($g(A)=1$) this contradicts $p\nmid g(A)$, therefore $\det A\neq4k^{3}$ by
Lemma V.

\subsection{Bound and tightness ($-$)}

\S 2.2 and \S 2.3 give $|\det A|\leq4k^{3}-k$. \ The symmetric extremizer and
its border are found via $(\star)$ with $d=0$, $\mathbf{v}=(1,1,0)$,
$\mathbf{u}=(1,1,-1)$, for which $\mathbf{v}^{\!\top}\operatorname{adj}%
(A_{k}^{-})\mathbf{u}=-1$:
\[
A_{k}^{-}=%
\begin{pmatrix}
k & k & k\\
k & -k & k-1\\
k & k-1 & -k
\end{pmatrix}
,\qquad M_{A}^{-}=%
\begin{pmatrix}
k & k & k & 1\\
k & -k & k-1 & 1\\
k & k-1 & -k & -1\\
1 & 1 & 0 & 0
\end{pmatrix}
.
\]
Here $\det A_{k}^{-}=4k^{3}-k$ by direct expansion and $\operatorname{adj}%
(A_{k}^{-})$ has a coprime pair of entries $2k-1$ and $k$, giving $g(A_{k}%
^{-})=1$; while $\Vert M_{A}^{-}\Vert=k$ and $\det M_{A}^{-}=1$ for all
$k\geq2$.

\section{Comparison\medskip}%

\[%
\begin{array}
[c]{c|cc}
& \mathbf{nonnegative}\{0,\dots,k\} & \mathbf{signed}\{-k,\dots,k\}\\\hline
\text{cube max} & 2k^{3} & 4k^{3}\\
\text{submatrix bound} & 2k^{3}-2k^{2}+k & 4k^{3}-k\\
\text{gap below cube max} & 2k^{2}-k & k\\
\text{achievability tool} & \text{Lemma C}\rightarrow k^{2}\mid\det
\;\;\;\;\; & \text{Lemma P + alignment}\\
\text{rigid survivors above bound} & 2k^{3}-k^{2},\ 2k^{3} & 4k^{3}\\
\text{killed by} & k\mid g(A) & p\mid g(A)\\
\text{extremal family} & A_{k}^{+},\ B_{k}^{+} & A_{k}^{-}%
\end{array}
\]

Both proofs are instances of "achievability narrows the above-bound spectrum
to rigid values, all carrying a prime factor of $k$ in their minors; Lemma E
($g=1$) deletes them". The shared embedding property is the hinge; the
divergent engines (the $k$-divisibility of Lemma C versus the minor-gap of
Lemma P) are exactly why two proofs, rather than one, give the clearest account.

\section{Closing Words}

Let $k\geq2$. Our main result is%
\[
\max_{M}\;\max_{A}\;\left\vert \det A\right\vert =\left\{
\begin{array}
[c]{ll}%
k(2k^{2}-2k+1) & \text{if }M\in\{0,1,\ldots,k\}^{4\times4},\\
k\left(  2k-1\right)  \left(  2k+1\right)  & \text{if }M\in\{-k,\ldots
,-1,0,1,\ldots,k\}^{4\times4}%
\end{array}
\right.
\]
where optimization is taken over all $3\times3$ submatrices $A$ of $4\times4$
unimodular matrices $M$ satisfying $\left\Vert M\right\Vert =k$. \ This has
implications for matrix inverses.

\textbf{Corollary.}%
\[
\max_{M}\;\left\Vert M^{-1}\right\Vert =\left\{
\begin{array}
[c]{ll}%
k(2k^{2}-2k+1) & \text{if }M\in\{0,1,\ldots,k\}^{4\times4},\\
k\left(  2k-1\right)  \left(  2k+1\right)  & \text{if }M\in\{-k,\ldots
,-1,0,1,\ldots,k\}^{4\times4}%
\end{array}
\right.
\]
where optimization is taken over all $4\times4$ unimodular matrices $M$
satisfying $\left\Vert M\right\Vert =k$. \ 

\textbf{Proof.} \ Because $\det M=\pm1$, the entries of $M^{-1}=\pm
\operatorname{adj}(M)$ are precisely the (signed) $3\times3$ minors of $M$,
thus $\|M^{-1}\|=\max_{A}|\det A|$ and the claim follows from the Theorem.
$\blacksquare$

For $3\times3$ matrices with entries in $\{0,1,\dots,k\}$, the largest
attainable determinant is the Hadamard-type ceiling $2k^{3}$, and Spectral
Rigidity shows that the only values surviving in the window $(\beta_{k}%
^{+},2k^{3}]$ are the two multiples of $k^{2}$, namely $2k^{3}$ and
$2k^{3}-k^{2}$. Below that rigid pair the spectrum becomes richer, and the
question "what are the largest \textit{genuine} (primitive, $g=1$)
determinants?" turns out to have a somewhat surprising answer. The first five
appear to be
\[
\beta_{k}^{+}=k(2k^{2}-2k+1)\ \succ\ 2k^{2}(k-1)\ \succ\ k^{3}+(k-1)^{3}%
\ \succ\ k(2k^{2}-3k+2)\ \succ\ (2k^{2}-k+1)(k-1)
\]
with successive gaps%
\[
k,\qquad k^{2}-3k+1,\qquad k-1,\qquad1.
\]
The determinants leave the $k$-sublattice already at the third step
($k^{3}+(k-1)^{3}\equiv-1\,(\operatorname{mod}k)$). It is produced not by the
inner-optimal rows whose large product-sum forces $k\mid\det$ -- the mechanism
governing the top of the multiple-of-$k$ ladder -- but instead by the rows of
a circulant
\[
\left(
\begin{array}
[c]{ccc}%
k & 0 & k-1\\
k-1 & k & 0\\
0 & k-1 & k
\end{array}
\right)
\]
whose determinant is simply the sum of two cubes. This level is exactly what
the original product-sum ladder, restricted to multiples of $k$, could not see.

These are \textit{exploratory results}: the constructions give rigorous lower
bounds (each value is exhibited by an explicit matrix), and the computations
confirm there is nothing between them up through $k=7$, but a proof that
$k^{3}+(k-1)^{3}$ is \textit{exactly} the third value for all $k\geq3$ -- and
a full characterization of the descent below it -- remains open. We would
warmly encourage others to extend the search to $k=8,9,\dots$, to test whether
the "first primitive level $=$ sum of two consecutive cubes" phenomenon
persists, and to chase down the analogous signed-alphabet ladder. There is a
great deal of structure here still waiting to be explained, and the chase is
an inviting one.

\section{Appendix:\ Base Case}

Lemma E specializes: for a $2\times2$ submatrix $A$ of a $3\times3$ unimodular
$M$, the gcd of the $1\times1$ minors -- i.e., of the four entries of $A$ --
equals $1$. (In the bordering identity $(\star)$, $\operatorname{adj}(A)$ now
has the \textit{entries} of $A$ as its components.)

\subsection{Nonnegative case, $A\in\{0,\dots,k\}^{2\times2}$}

The cube max is $k^{2}$, and it is already attainable with coprime entries:
\[
A=%
\begin{pmatrix}
k & 1\\
0 & k
\end{pmatrix}
,\quad\det A=k^{2},\quad\gcd=1,\qquad M=%
\begin{pmatrix}
k & 1 & 0\\
0 & k & 1\\
1 & 0 & 0
\end{pmatrix}
.
\]
Thus no correction occurs: the bound is exactly the cube max $k^{2}$.
(Contrast with $n=4$, where the maximizing $\{0,k\}$-configuration is forced
to have $g=k$, producing the $-2k^{2}+k$ deficit. At $n=3$ the maximizer
escapes the obstruction.)

\subsection{Signed case, $A\in\{-k,\dots,k\}^{2\times2}$}

Now $\det A$ \textit{is} a $2\times2$ minor, so Lemma P applies directly: no
value in $(2k^{2}-k,\,2k^{2})$ is attained. The cube max $2k^{2}$ forces all
four entries $=\pm k$, hence $\gcd=k>1$, excluded by Lemma E. Therefore
\[
|\det A|\leq2k^{2}-k=k(2k-1),
\]
attained with coprime entries:
\[
A=%
\begin{pmatrix}
k & k\\
-(k-1) & k
\end{pmatrix}
,\quad\det A=2k^{2}-k,\quad\gcd=1,\qquad M=%
\begin{pmatrix}
k & k & -1\\
-(k-1) & k & 1\\
0 & 1 & 0
\end{pmatrix}
.
\]
This is the whole proof -- a single application each of Lemma P and Lemma E,
with no deficit count and no induction. It is the most transparent possible
illustration of the paper's architecture. \ Of course we have, for $k\geq2$,%
\[
\max_{M}\;\left\Vert M^{-1}\right\Vert =\left\{
\begin{array}
[c]{ll}%
k^{2} & \text{if }M\in\{0,1,\ldots,k\}^{3\times3},\\
k\left(  2k-1\right)  & \text{if }M\in\{-k,\ldots,-1,0,1,\ldots,k\}^{3\times3}%
\end{array}
\right.
\]
where optimization is taken over all $3\times3$ unimodular matrices $M$
satisfying $\left\Vert M\right\Vert =k$. \ 

Ternary counts are given in \cite{S3-matrx3, S4-matrx3} but the subject is
undeveloped, as far as is known.

Our extremal $3\times3$ and$\ 4\times4$ matrices may be viewed as
fixed-dimension unimodular analogs of the anti-Hadamard constructions of
Graham \&\ Sloane \cite{GS-matrx3} and Alon \&\ V\~{u} \cite{AV-matrx3}. \ 

An early draft of our paper employed Gram matrices and the work of
R\'{o}\.{z}a\'{n}ski \&\ Witu\l a \cite{RW-matrx3}. \ 

Extremal matrices first appeared in \cite{Fn-matrx3}, as an afterthought,
under tighter restrictions (none of the entries in $M$ and none of the entries
in $M^{-1}$ were allowed to be zero). \ It seemed important to restate our
conjecture in full generality, removing the zerofree condition. \ Also, the
case $k=1$ was empty in \cite{Fn-matrx3} but it is meaningful here.

\section{Acknowledgements}

Microsoft Copilot assisted in writing CUDA programs, which opened the door for
me to run intensive experimental jobs on Nvidia GPUs. \ Anthropic Claude Opus
4.8 played a crucial role in turning my conjectures into theorems. \ The
creators of Mathematica, as well as administrators of the MIT Engaging
Cluster, earn my gratitude every day.\pagebreak

\section{Addendum I}

Returning to the first paragraph, we exhibit matrices indicating that
$\gamma_{8}^{-}(1)\geq442$ and $\gamma_{10}^{+}(1)\geq110$:%
\[%
\begin{array}
[c]{ccc}%
\left(
\begin{array}
[c]{cccccccc}%
0 & -1 & 1 & -1 & 0 & 1 & 1 & 1\\
-1 & 1 & 0 & 1 & -1 & 1 & 1 & 0\\
0 & 1 & 1 & -1 & -1 & -1 & 1 & 1\\
-1 & 0 & 1 & -1 & -1 & 1 & 0 & 1\\
1 & 1 & 1 & -1 & -1 & 1 & 1 & -1\\
-1 & 1 & 0 & -1 & 1 & 1 & -1 & 1\\
1 & 0 & 1 & 1 & -1 & 1 & -1 & 1\\
-1 & -1 & 0 & -1 & -1 & 0 & -1 & -1
\end{array}
\right)  , &  & \left(
\begin{array}
[c]{cccccccccc}%
1 & 0 & 1 & 0 & 1 & 1 & 0 & 1 & 1 & 0\\
1 & 0 & 1 & 1 & 0 & 0 & 1 & 1 & 0 & 0\\
0 & 1 & 1 & 1 & 0 & 0 & 0 & 1 & 1 & 1\\
0 & 1 & 0 & 1 & 0 & 0 & 1 & 0 & 1 & 1\\
0 & 1 & 0 & 1 & 0 & 1 & 0 & 1 & 1 & 0\\
1 & 1 & 0 & 0 & 0 & 0 & 1 & 0 & 1 & 0\\
1 & 0 & 0 & 1 & 1 & 1 & 0 & 0 & 1 & 1\\
1 & 1 & 0 & 0 & 0 & 1 & 0 & 1 & 0 & 1\\
0 & 0 & 0 & 0 & 0 & 0 & 1 & 1 & 1 & 1\\
0 & 1 & 1 & 0 & 1 & 1 & 1 & 0 & 0 & 1
\end{array}
\right)  .
\end{array}
\]
Data on $\gamma_{n}^{-}(k)$ appears in \cite{Fn-matrx3}, with an additional
zerofree constraint. \ We shall therefore focus on $\gamma_{n}^{+}(k)$. \ The
first four terms of the sequence $\gamma_{3}^{+}(2)=4$, $\gamma_{4}^{+}%
(2)=10$, $\gamma_{5}^{+}(2)=30$, $\gamma_{6}^{+}(2)=130$, $\gamma_{7}%
^{+}(2)=454$, ... match the terms $\gamma_{4}^{-}(1)$, $\gamma_{5}^{-}(1)$,
$\gamma_{6}^{-}(1)$, $\gamma_{7}^{-}(1)$ \textit{exactly}, a striking
concidence across different alphabets and different dimensions. \ Might
serendipity collapse here, i.e., could it be true that $\gamma_{8}^{-}%
(1)<454$?\medskip

Further $\gamma_{5}^{+}(k)$, $\gamma_{5}^{-}(k)$, $\gamma_{6}^{-}(k)$
calculations for $k\geq2$ lead to the following.\medskip

\textbf{Conjecture. \ }We have
\[
\gamma_{5}^{+}(k)=\left\{
\begin{array}
[c]{lll}%
k(3k^{3}-3k^{2}+2k-1) &  & \text{if }k\equiv2\,\left(  \operatorname{mod}%
3\right)  ,\\
3k^{2}(k^{2}-k+1) &  & \text{otherwise;}%
\end{array}
\right.
\]%
\[
\gamma_{5}^{-}(k)=k(4k-1)(4k^{2}-2k+1);\medskip
\]%
\[
\gamma_{6}^{-}(k)=\left\{
\begin{array}
[c]{lll}%
k(4k^{2}-2k+1)(12k^{2}-2k-1) &  & \text{if }k\equiv1\,\left(
\operatorname{mod}5\right)  ,\\
k(2k-1)(6k-1)(4k^{2}+1) &  & \text{otherwise.}%
\end{array}
\right.
\]
A formula for $\gamma_{6}^{+}(k)$ will be possible only if existing data%
\[
\gamma_{6}^{+}(k)=\left\{
\begin{array}
[c]{ccc}%
30 &  & \text{if }k=2,\\
1011 &  & \text{if }k=3,\\
4420 &  & \text{if }k=4,\\
13705 &  & \text{if }k=5,\\
34854 &  & \text{if }k=6,\\
76055 &  & \text{if }k=7,\\
150280 &  & \text{if }k=8\text{ }%
\end{array}
\right.
\]
can be suitably extended. \ All that is presently known: $\gamma_{6}%
^{+}(9)\geq190332$ and $\gamma_{6}^{+}(10)\geq263610$, which correspond
to\pagebreak\
\[%
\begin{array}
[c]{ccc}%
\begin{array}
[c]{c}%
\left(
\begin{array}
[c]{cccccc}%
9 & 9 & 9 & 9 & 6 & 1\\
9 & 9 & 1 & 0 & 9 & 0\\
0 & 8 & 0 & 9 & 9 & 0\\
9 & 0 & 0 & 9 & 8 & 0\\
0 & 0 & 9 & 1 & 9 & 0\\
8 & 0 & 1 & 6 & 7 & 0
\end{array}
\right)  ,
\end{array}
&  & \left(
\begin{array}
[c]{cccccc}%
10 & 10 & 10 & 10 & 10 & 1\\
10 & 1 & 0 & 0 & 10 & 1\\
10 & 0 & 0 & 10 & 1 & 0\\
10 & 0 & 10 & 1 & 0 & 0\\
10 & 10 & 1 & 0 & 0 & 0\\
1 & 7 & 4 & 5 & 0 & 0
\end{array}
\right)  .
\end{array}
\]

\section{Addendum II}

In fact, $\gamma_{8}^{-}(1)=454$. \ The matrix
\[
\left(
\begin{array}
[c]{cccccccc}%
1 & 1 & 1 & 1 & 1 & 1 & 1 & 1\\
1 & 1 & -1 & -1 & 1 & 1 & -1 & 0\\
0 & -1 & 1 & -1 & 1 & 1 & -1 & -1\\
1 & 1 & 1 & -1 & -1 & -1 & 0 & 1\\
1 & -1 & -1 & -1 & 0 & -1 & 1 & 0\\
-1 & 1 & 0 & -1 & 1 & -1 & 1 & 0\\
-1 & 0 & -1 & -1 & -1 & 1 & 1 & 0\\
-1 & -1 & -1 & -1 & -1 & 0 & -1 & -1
\end{array}
\right)
\]
indicates that $\gamma_{8}^{-}(1)\geq454$, and our computer-assisted
certification establishes that no larger value is possible. \ We can build
nearly orthogonal integer rows (akin to $\pm1$ Hadamard-like patterns), close
to saturating the bound used in our code by the pruner. \ By a
\textbf{borderable core} is meant a submatrix $A$ that admits a border -- a
row, a column, and a corner entry -- completing it to a unimodular matrix $M$.
\ If a determinant-preserving bijection between borderable $\{0,1,2\}$ cores
of order $n-1$ and borderable $\{-1,0,1\}$ cores of order $n$ exists, then
$\gamma_{n}^{+}(2)=\gamma_{n+1}^{-}(1)$ would be true.

Likewise,\ $\gamma_{10}^{+}(1)=120$. \ The matrix
\[
\left(
\begin{array}
[c]{cccccccccc}%
1 & 1 & 1 & 1 & 1 & 1 & 0 & 0 & 0 & 1\\
1 & 1 & 0 & 1 & 0 & 0 & 0 & 1 & 1 & 0\\
0 & 0 & 0 & 1 & 1 & 1 & 1 & 0 & 1 & 1\\
0 & 1 & 1 & 0 & 0 & 1 & 1 & 0 & 1 & 0\\
1 & 1 & 0 & 0 & 1 & 0 & 1 & 0 & 1 & 0\\
1 & 0 & 1 & 1 & 0 & 0 & 1 & 0 & 1 & 0\\
0 & 1 & 0 & 1 & 1 & 0 & 1 & 1 & 0 & 0\\
0 & 0 & 1 & 0 & 1 & 0 & 0 & 1 & 1 & 0\\
1 & 0 & 0 & 0 & 0 & 1 & 1 & 1 & 0 & 0\\
1 & 1 & 1 & 0 & 0 & 1 & 1 & 0 & 0 & 0
\end{array}
\right)
\]
and certification play roles here. \ Two rows can never be \textquotedblleft
orthogonal with cancellation\textquotedblright\ the way that signed rows can;
thus the pruning bound on which our code depends is far from tight. \ Such
looseness inescapably leads to prolonged runtimes.

\section{Addendum III}

We conjecture that, if $2\leq k\not \equiv 6$ $\left(  \operatorname{mod}%
7\right)  $, then%
\[
\gamma_{7}^{-}(k)=k(4k^{2}+1)(40k^{3}-20k^{2}+6k-1).
\]
An allied expression for $k\equiv6\,\left(  \operatorname{mod}7\right)  $
remains open because, while $\gamma_{7}^{-}(6)=6872190$, the value of
$\gamma_{7}^{-}(13)$ is currently inaccessible.

A formula for $\gamma_{7}^{+}(k)$ will be possible only if existing data%
\[%
\begin{array}
[c]{ccc}%
\gamma_{7}^{+}(2)=454=\gamma_{8}^{-}(1), &  & \gamma_{7}^{+}(3)=5760
\end{array}
\]
can be suitably extended. \ All that is presently known: $\gamma_{7}%
^{+}(4)\geq21662$ and $\gamma_{7}^{+}(5)\geq72831$, which correspond to%
\[%
\begin{array}
[c]{ccc}%
\left(
\begin{array}
[c]{ccccccc}%
1 & 0 & 1 & 4 & 4 & 4 & 4\\
0 & 4 & 3 & 0 & 4 & 3 & 0\\
4 & 3 & 1 & 0 & 0 & 4 & 4\\
4 & 3 & 1 & 4 & 2 & 1 & 0\\
0 & 4 & 2 & 3 & 1 & 0 & 4\\
4 & 2 & 3 & 0 & 4 & 0 & 1\\
3 & 3 & 0 & 3 & 0 & 2 & 0
\end{array}
\right)  , &  & \left(
\begin{array}
[c]{ccccccc}%
4 & 2 & 2 & 0 & 2 & 1 & 3\\
0 & 5 & 5 & 4 & 1 & 2 & 3\\
2 & 0 & 5 & 4 & 5 & 0 & 0\\
4 & 5 & 5 & 3 & 0 & 4 & 0\\
0 & 5 & 1 & 5 & 5 & 5 & 0\\
5 & 0 & 5 & 0 & 1 & 5 & 5\\
5 & 4 & 0 & 0 & 5 & 0 & 4
\end{array}
\right)  .
\end{array}
\]

Finally, we have%
\[%
\begin{array}
[c]{ccc}%
\gamma_{8}^{-}(2)=62550, &  & \gamma_{8}^{-}(3)=1129299
\end{array}
\]
and it appears that -- barring a significant theoretical breakthrough --
inferring $\gamma_{8}^{-}(k)$ for arbitrary $k$ by computation alone will be formidable.

\end{document}